\documentclass[
11pt,%
tightenlines,%
twoside,%
onecolumn,%
nofloats,%
nobibnotes,%
nofootinbib,%
superscriptaddress,%
noshowpacs,%
centertags]%
{revtex4-2}
\usepackage{ljm}
\usepackage{color}

\usepackage{stmaryrd}

\newcommand{\CC}{\ensuremath{{\mathbb C}}}
\newcommand{\CW}{\ensuremath{{\widehat{\mathbb C}}}}
\newcommand{\RR}{\ensuremath{{\mathbb R}}}
\newcommand{\ZZ}{\ensuremath{{\mathbb Z}}}
\newcommand{\NN}{\ensuremath{{\mathbb N}}}
\newcommand{\XX}{\ensuremath{{\mathbb X}}}
\renewcommand{\Re}[1]{{\mathfrak{Re}\left(#1\right)}}
\renewcommand{\Im}[1]{{\mathfrak{Im}\left(#1\right)}}

\newcommand{\HH}{\ensuremath{{\mathbb H}}}
\font\myfont=cmr10 at 12pt  \newcommand{\e}{{\text{\myfont e}}}
\newcommand{\del}[2]{\frac{\partial #1}{\partial #2}}
\newcommand{\ent}[2]{{}_{#1}\mathscr{E}_{#2}}
\renewcommand{\SS}{\mathbb{S}}

\usepackage{graphicx}
\makeatletter
\newcommand*\bigcdot{\mathpalette\bigcdot@{.5}}
\newcommand*\bigcdot@[2]{\mathbin{\vcenter{\hbox{\scalebox{#2}{$\m@th#1\bullet$}}}}}
\makeatother

\usepackage{enumitem}
\setlist[enumerate]{
wide, 
nosep, 
labelwidth=*,
labelindent=0cm,
leftmargin=0cm}

\begin{document}
\titlerunning{Integrability and complex structures} 
\authorrunning{Le\'on--Gil, Muci\~no--Raymundo} 

\title{Integrability and adapted complex 
\\ structures to smooth  vector fields on the plane}

\author{\firstname{Gaspar} \surname{Le\'on--Gil}}
\email[E-mail: ]{leon.gil.gaspar@gmail.com} \affiliation{
Tecnol\'ogico de Tac\'ambaro, Tac\'ambaro, Michoac\'an 61650, M\'exico
}

\author{\firstname{Jes\'us} \surname{Muci\~no--Raymundo}}
\email[E-mail: ]{muciray@matmor.unam.mx} \affiliation{Centro de Ciencias Matem\'aticas, UNAM, Morelia, Michoac\'an 58089, M\'exico}



 \received{March 11, 2020; revised March 18, 2020; accepted April 10, 2020}

\begin{abstract} 
Singular complex analytic vector fields on the Riemann surfaces enjoy several
geometric properties (singular means that poles and essential singularities 
are admissible). We describe relations between 
singular complex analytic vector fields $\XX$ and smooth vector fields $X$.  
Our approximation route studies three 
integrability notions for real smooth vector fields $X$
with singularities on the plane or the sphere.
The first notion is related to Cauchy--Riemann 
equations, we say that
a vector field $X$ admits an adapted complex structure 
$J$ if there exists a singular complex analytic 
vector field $\XX$ on the plane provided with this complex structure,
such that $X$ is the real part of $\XX$.
The second integrability notion for $X$
is the existence of a 
first integral $f$, smooth and having non vanishing differential 
outside of the singularities of $X$. 
A third concept is that $X$ admits a  
global flow box map 
outside of its singularities,  {\it i.e.} the vector field $X$
is a lift of the trivial horizontal vector field, under a diffeomorphism. 
We study the relation between the three notions.
Topological obstructions (local and global) to the three
integrability notions are described. 
A construction of singular complex analytic vector fields 
$\XX$ using canonical invariant regions is provided. 
\end{abstract}

\maketitle


\section{Introduction}
\label{intro}

Our aim is to characterize dynamically the 
real vector fields that coincide with the real 
parts of singular complex analytic vector fields.
Let 
$\SS^2 = \RR^2\cup \{ \infty \} $,
we consider a set $\mathfrak{S} \subset \RR^2$, 
having at most a finite number of accumulation points.  
Let $X \in \mathfrak{X}^\infty (\RR^2  \backslash \mathfrak{S})$ 
be a $C^\infty$ vector field with two kind of singularities

\noindent
$\bigcdot$ smooth zeros at $\mathcal{Z}(X) \subset 
\RR^2 \backslash \mathfrak{S}$, and 

\noindent 
$\bigcdot$ 
wide singularities, {\it i.e.} points in $\mathfrak{S}$ where   
$X$ is undefined or non $C^\infty$.

\noindent
Consider $\mathfrak{P}= \mathfrak{S} \cup \mathcal{Z}(X) $,
thus 
$\RR^2 \backslash \mathfrak{P}$ is a plane with topological punctures.
There are plenty of complex structures $\{ J \}$ such that 
$(\RR^2 \backslash \mathfrak{P}, J)$ is a Riemann surface.

\medskip 

\noindent 
{\it
Let 
$X\in \mathfrak{X}^\infty (\RR^2 \backslash \mathfrak{S})$,
under which conditions
are there a complex structure $J$ and 
a singular complex analytic vector field $\XX$
on the Riemann surface $(\RR^2 \backslash  \mathfrak{P}, J)$, 
such that }
\begin{equation}\label{ecuacion-de-estructura-compleja-adaptada}
\rho X= \Re{\XX} \ ?
\end{equation}

\smallskip

\noindent 
Here $\rho$ is a  $C^\infty$
non vanishing reparametrization on $\RR^2 \backslash \mathfrak{P}$
and
$\Re{\XX}$ denotes the real part of  $\XX$. The adjective {\it singular} 
means that $\XX$ may have poles and  essential singularities 
at $\mathcal{S} \cup \infty$.  
For affirmative cases, we say that $J$ is an  
{\it adapted complex structure} to the vector field $X$.
In fortunate situations, $J$ defines conformal punctures 
at $\mathfrak{P}$ and a maximal
Riemann surface $(\RR^2, J)$ emerges (conformally equivalent to the 
plane $\CC$ or the Poincar\'e disk $\Delta$).
The problem 
\eqref{ecuacion-de-estructura-compleja-adaptada}  
for meromorphic vector fields $\XX$ on compact orientable surfaces, 
was studied in \cite{Mucino}.

\smallskip 

Secondly, we define that $X$ is {\it integrable}  if there exists 
an integrating factor $\mu$ 
and a Hamiltonian vector field $X_f$ 
of a function $f$, with  $df$ non vanishing on 
$\RR^2 \backslash \mathfrak{P}$,
such that 
\begin{equation}\label{ecuacion-de-integrabilidad}
\mu X = X_{f},
\end{equation}
\noindent 
here $\mu$ and $f$ are $C^\infty$ on 
$\RR^2 \backslash \mathfrak{P} $.
{\it Under which conditions $X$ is integrable?}

\smallskip 

\noindent 
In fact, the 
non vanishing of $df$ on 
$\RR^2 \backslash \mathfrak{P}$
is the innovative condition, looking to other integrability concepts
in the literature.

Our third notion is as follows.
A vector field 
$X$ {\it admits a global flow box} if there exists 
an scaling factor $\rho$
and a probably multivalued map
such that
\begin{equation}\label{ecuacion-de-caja-de-flujo-global}
(g,f): \RR^2 \backslash \mathfrak{P} 
\longrightarrow \RR^2,  
\ \
(g,f)_* (\rho X) = \del{}{\tau},
\end{equation}

\noindent
here $\rho$ and $(g,f)$ are $C^\infty$ on 
$\RR^2 \backslash \mathfrak{P}$.
{\it Under which conditions $X$ admits a global flow box?}

\smallskip

\noindent 
Geometrically, \eqref{ecuacion-de-caja-de-flujo-global} 
means that, outside of its singularities the vector field $X$
is {\it a lift } of the horizontal vector field
$\partial / \partial \tau$ on $\RR^2$ under a probably multivalued map. 
Liftable vector fields appear in     
singularity theory V.\,I.\,Arnold \cite{Arnold1}  p.\,561 or 
A.\,A. du Plessis {\it et al.} \cite{Plessis-Wall} p.\,120, and
in Riemann surface theory \cite{Boothby2}, \cite{Bott}.

\medskip

In our framework, a
{\it singular complex analytic, 
additively automorphic, single valued or multivalued
function}  
$\Psi : 
(\RR^2 \backslash \mathfrak{P}, J)  \longrightarrow \CC$ 
has meromorphic local branchs  
and single valued differential $d\Psi$.
In some cases, $\Psi$ extends to $\mathfrak{S} \cup \infty$ meromorphically 
and/or 
having  essential singularities. 
A source of nice affirmative examples for equations 
\eqref{ecuacion-de-estructura-compleja-adaptada}--\eqref{ecuacion-de-caja-de-flujo-global}
is as follows.

\begin{theorem}[A dictionary]
\label{teorema-1-introduccion-diccionario} 
There exists a one to one correspondence 
between: 

\begin{enumerate}[label=\arabic*)]
\item 
Singular complex analytic vector fields $\XX$ on 
a Riemann surface 
$(\RR^2 \backslash \mathfrak{P}, J)$.

\item 
Singular complex analytic 
(additively automorphic, single valued or multivalued)
functions 

\centerline{$\Psi = g + \sqrt{-1}f + c$ \ on 
$(\RR^2 \backslash \mathfrak{P}, J)$, $c \in \CC$.}

\item
Integrable $C^\infty$ real vector fields $X$ on 
$\RR^2 \backslash \mathfrak{P}$.

\item 
$C^\infty$ real vector fields $X$ admitting 
a global flow box 
$(g,f)$ on  $\RR^2 \backslash \mathfrak{P}$.
\end{enumerate}

\smallskip 

\noindent 
Moreover, the correspondence is such that 
$$
\rho X=\Re{\XX }, 
\ \ 
\Psi_* \XX = \del{}{t}, 
\ \   
f= \Im{\Psi} 
\ \hbox{ and } \  
df(X)=0,  
$$
here $t= \tau+ \sqrt{-1} \sigma$ denotes the complex time of $\XX$.
\end{theorem}

What does a singular complex analytic vector field 
$\XX$ look like?
As usual let us denote, 
$\mathbb{H}^2$ the open half plane,
$\Delta_R$ the open disk of radius $R\geq 1$ in $\CC$, and 
$\overline{(\ \, )}$ the topological closure.
A dynamical/constructive 
characterization of vector fields
$\XX$ is as follows.
 
\begin{theorem}[Decomposition for singular complex analytic vector fields] 
\label{teorema-2-introduccion-Fatou-decomposition}
Let $M \subseteq \SS^2$ be an open connected surface. 
\begin{enumerate}[label=\arabic*)]
\item 
Assume that $M$ is
obtained by the paste of a finite or infinite number of closed canonical regions of type
$$
\left( 
\overline{\mathbb{H}}^2, \del{}{z} \right)
,\
\left( \{ 0 \leq \Im{z} \leq h \}, \frac{\partial}{\partial z} \right)
,\
\left( \overline{\Delta}_1, \frac{-2 \pi iz}{r} \del{}{z}\right)
,\
\left( \overline{\Delta}_R \backslash \Delta_1, 
\frac{-2 \pi iz}{r} \del{}{z} \right), 
$$

\noindent 
$h,\, r, R \in \RR^*$. 
Then there exist a complex structure $J$ and 
a singular complex analytic vector field $\XX$ on $(M, J)$, 
extending the vector fields of the canonical regions.

\item 
Conversely, 
let $\XX$ be a singular complex analytic vector field on 
$(M, J)$, having
at most a locally finite set of real incomplete 
trajectories
$\{ z_\vartheta (\tau) \} $.

\noindent
Then $\XX$ admits a locally finite 
decomposition in  regions as above. 
\end{enumerate}
\end{theorem}

As a novel aspect
we present decompositions with an infinite number
of pieces. 
In order to study the questions
\eqref{ecuacion-de-integrabilidad} and \eqref{ecuacion-de-caja-de-flujo-global}
we require  some concepts.

\noindent 
A vector field 
$X \in \mathfrak{X}^\infty (\RR^2 \backslash \mathfrak{S})$ 
has the following remarkable trajectory sets:

\vskip0.1cm 

\noindent
The {\it separatrix trajectories 
$\Gamma(X)=\{q_j\} \cup \{ \zeta_\varsigma \} $ of $X$} are

\noindent $\bigcdot$
points\footnote{Here we abuse of the 
notation, since a non smooth point 
is not a trajectory of $X$.} 
$q_j \in (\mathfrak{P} \cup \{\infty \})
\backslash \{\hbox{topological centers, sources or sinks}\} $,

\noindent 
$\bigcdot$
non stationary $\zeta_\varsigma=\zeta_\varsigma(\tau)$ trajectories, such that 
do not admit an open neighbourhood in $\RR^2 \backslash \mathfrak{P}$
filled by 
trajectories having the same topology.

\noindent 
The {\it separatrix skeleton} $X$
is $\Gamma (X)= \{ q_j \} \cup \{ \zeta_\varsigma \}$, 
a graph with 
vertices $\{ q_j \}$ and edges $\{ \zeta_\varsigma \}$. 

\smallskip

\noindent 
The {\it attractors $\mathcal{A}(X)$ of $X$} are 

\noindent 
$\bigcdot$
points $q_j \in \mathfrak{P}$ which 
admits topological parabolic sectors
(in particular, topological sources or sinks),

\noindent 
$\bigcdot$
periodic trajectories and polycycles (union of separatrices 
$\zeta_\varsigma$, points $q_j \in \mathfrak{P} \cup 
\{\infty \}$  homeomorphic to a 
circle $\SS^1 \subset \RR^2$),  
whose holonomy germ is different from the identity.

\smallskip

\begin{corollary}[Global flow box for 
$C^\infty$ real vector fields]
\label{teorema-3-introduccion-caracterizacion-GFB} 
Let $\XX$ be a complex analytic vector field on 
$(\RR^2 \backslash \mathfrak{P}, J)$
with a locally finite set of incomplete real trajectories. 
Then 
the  vector field  
$X=\Re{\XX} \in \mathfrak{X}^\infty (\RR^2 \backslash \mathfrak{S} )$ 
be a vector field that satisfies the following conditions:

\begin{enumerate} 
\item 
The separatrices $\Gamma(X)$ 
determine a locally finite set 
of trajectories
in $\RR^2 \backslash \mathfrak{P}$.

\item 
The periodic trajectories and polycycles in $\Gamma(X)$
have $C^\infty$ identity as holonomy or first return 
maps. 

\item
The holonomy of each hyperbolic sector
at $q_j \in \mathfrak{P}  \cup \{ \infty\}$
is a $C^\infty $ diffeomorphism.

\item
For each $q_j \in \mathfrak{P}\cup \{ \infty \}$ 
a topological multi--saddle 
with $2k + 2 \geq 2 $ topological hyperbolic sectors, 
$X$ is $C^\infty$ equivalent to 

\centerline{$
\Re { \frac{1}{z^k} \del{}{z}} = 
\Re{\frac{1}{z^k}} \del{}{x} + \Im {\frac{1}{z^ k}} \del{}{y}, \ \ \ 
k\in \NN \cup \{ 0\},
$}

\noindent 
in a punctured neighbourhood of $q_j$.
\end{enumerate}
\end{corollary}

In particular, the limit cycles 
are obstructions in order to get affirmative 
solutions questions
\eqref{ecuacion-de-estructura-compleja-adaptada}--\eqref{ecuacion-de-caja-de-flujo-global}.
About our hypothesis 
``$\XX$ having locally finite set of incomplete real trajectories'':
on $\CW$ this set is finite if and only if $X$ is rational, 
see Corollary \ref{caracterizacion-racional}.
In particular the existence of an essential singularity for $\XX$
implies an infinite number of incomplete real trajectories.  
Let us recall that our vector fields enjoy  
two geometric properties.

\begin{corollary}
\label{corolario-1-consecuencias-geometricas} 
Let 
$X \in \mathfrak{X}^\infty (\RR^2 \backslash \mathfrak{S})$ be a vector field as 
in Theorem 
\ref{teorema-3-introduccion-caracterizacion-GFB}.

\begin{enumerate}

\item
There exists a $C^\infty $ flat Riemannian metric
$g_X$ on $\RR^2 \backslash 
\big( \mathcal{Z}(X) \cup \mathfrak{S} \cup \mathcal{A} (X) \big)$, 
such that $X$ is a geodesic vector field. 

\item
$X$ is one of the two linearly independent infinitesimal generators 
of a $C^\infty $ local $(\RR^2, +)$--action on 
$\RR^2 \backslash \big( \mathcal{Z}(X) \cup \mathfrak{S} \cup \mathcal{A} (X) \big)$.
\end{enumerate}
\end{corollary}

Among the families of $C^\infty $ vector fields satisfying 
the hypothesis in 
Theorem 
\ref{teorema-3-introduccion-caracterizacion-GFB}
there are; 
Hamiltonian vector fields $X_f$ and 
gradient $C^\infty$ vector 
fields $\nabla f$, from $f \in C^\infty(\RR^2, \RR)$, 
having all its zeros of Morse type.

The Uniformization Theorem asserts that any complex structure $J$ 
on $\RR^2$
makes it conformally equivalent to $\CC$ or the Poincar\'e
disk $\Delta$; however the recognizing problem is hard. 
Using vector fields $X$ with adapted complex structures
$(\RR^2, J)$ as in Theorem \ref{teorema-1-introduccion-diccionario},  
we want to recognize the induced complex analytic structure. 
Let us define that $X$ has a {\it finite trajectory gap} 
if in $\big( (\RR^2, J), X \big)$ there exists a holomorphic local flow box
$\Psi: \mathcal{U} \subset \RR^2 \to   \CC$ 
such that the image of the ideal boundary of $\RR^2$ under 
$\Psi$ is a simple path $\beta$ in $\CC$, see
Definition \ref{gap} and Figure \ref{3-fig}.

\begin{corollary}\label{corolario-2-introduccion-tipo-conforme}
Let $X \in \mathfrak{X}^\infty (\RR^2 \backslash \mathfrak{P})$ be real a vector field which is the real part of 
a singular complex analytic  vector field 
$\XX$ on $(\RR^2 \backslash \mathfrak{P}, J)$. 

\noindent 
1) If $X$ has a finite trajectory gap at a point $q$ of $\mathfrak{P}$, 
then $q$ is a conformal hole. 

\noindent 
2) If $J$ extends to $\mathfrak{P}$ ({\it i.e.} all the points in $\mathfrak{P}$ are conformal
punctures) and the respective 
$(\RR^2, J)$ has a finite trajectory gap, 
then it is biholomorphic to the Poincar\'e disk $\Delta$.
\end{corollary}

\noindent 
Convention. By notational simplicity, 
we work in the $C^\infty$ category, however all the
results remain valid under $C^1$ hypothesis.

The authors are very grateful with Alvaro Alvarez--Parrilla by
several illustrative conversations.

\section{First integrals and integrating factors}

We provide an explanation/review for the integrability equation
\eqref{ecuacion-de-integrabilidad}. Let 
$$
X(x,y) = a(x,y)\del{}{x} + b(x,y) \del{}{y} \ \in \ \mathfrak{X}^\infty (\RR^2 \backslash \mathfrak{S}),
$$
be a vector field, it has three associated objects:

\noindent 
$\bigcdot$
The sheaf of rings of first integrals of $X$  

\centerline{
$\mathcal{FI}(X) = \{ f \ \vert \  \mathcal{L}_X(f) \equiv 0 \},$
}

\noindent
under addition and multiplication. 
As a matter of record:
let $\{ U_{\tt j} \}$ be the cover by open sets
of $\RR^2 \backslash \mathfrak{S}$, 
the $C^\infty$
{\it sheaf of rings} 
$\mathcal{FI}(X)$ 
associates to each open $U_{\tt j}$ the 
ring  of the  $C^\infty$ first integrals
$f_\alpha : U_{\tt j}  \longrightarrow
\RR$ of $X$,
considering addition 
$f_\alpha + f_\beta$ and multiplication  $f_\alpha f_\beta$
as ring operations. 
Analogous $C^\infty$ sheaf notions 
apply for groups and 
Lie algebras below.

\noindent 
$\bigcdot$
The sheaf of groups of integrating factors of $X$ 

\centerline{
$\mathcal{IF}(X)= \{ \mu \ \vert \ \mu X \hbox{ is Hamiltonian}\}, $
}

\noindent
under addition $\mu_1 + \mu_2$.

\noindent
$\bigcdot$ 
The sheaf of Lie algebras of infinitesimal symmetries of $X$ 

\centerline{
$Sym(X) = \{ Y \ \vert \ [X, Y] = \nu X , \ \nu  \hbox{ function} \},$
}

\noindent
under the Lie bracket operation $[Y_1 ,Y_2]$.

\noindent 
The classical relations between these objects are described by the 
operators

\begin{equation}\label{trilogia-operadores}
\begin{picture}(150,47)
\put(65,0){$Sym (X)$ \ \ .}
\put(5,40){$\mathcal{FI}(X)$}
\put(127,40){$\mathcal{IF}(X)$}

\put(43,43){\vector(1,0){79}}
\put(16,33){\vector(3,-2){42}}
\put(152,33){\vector(-3,-2){44}}
\put(67,10){\vector(-3,2){38}}
\put(101,10){\vector(3,2){38}}

\put(77,47){$\mathfrak{c}_1 $}
\put(132,13){$\mathfrak{c}_5$}
\put(107,23){$\mathfrak{c}_4$}
\put(30,13){$\mathfrak{c}_2$}
\put(50,23){$\mathfrak{c}_3$}
\end{picture}
\end{equation}

\noindent
Here, the first operator and its inverse are 
\begin{equation}
\begin{array}{rcl}
\mathfrak{c}_1 :  \mathcal{FI}(X) & \longrightarrow  & \mathcal{IF}(X) 
\\
f & \longmapsto &\frac{-1}{a}\del{f}{y} = \frac{1}{b} \del{f}{x}
\\
\int^{(x, y)} \mu(- b dx + a dy) & \longmapsfrom &\mu. 
\end{array}
\end{equation}

\noindent 
The second operator is not canonical, we use
\begin{equation}\label{c-2}
\begin{array}{ccccccc}
\mathfrak{c}_2 : & \mathcal{FI}(X) & \longrightarrow  & Sym (X), \ \ \ 
& f & \longmapsto & 
Y= \frac{1}{-b \del{f}{x}+ a\del{f}{y}}
\left( -b \del{}{x} + a \del{}{y}\right)=
\frac{ \nabla f }{ \Vert  \nabla f \Vert ^2}  . \\
\end{array}
\end{equation}
The right equality follows from 
$\mathcal{L}_X= af_x + bf_y=0$,  
$a= -b (f_y/f_x)$, $b= -a(f_x / f_y)$ 
and a direct substitution in the first 
expression of $Y$. Note that
$\mathfrak{c}_2$ is not onto, since the resulting infinitesimal symmetries
$\mathfrak{c}_2 (f) =Y$ and $X$ are always orthogonal and
this condition is not fulfilled for every $Y \in Sym (X)$. 
In the reverse direction, 
an usual choice is  
\begin{equation}
\begin{array}{ccccccc}
\mathfrak{c}_3 : & Sym (X) & \longrightarrow  & \mathcal{FI}(X), \ \ \
& Y = c \del{}{x} + d\del{}{y} & \longmapsto & 
f (x,y) = \int^{(x,y)} \frac{-bdx + ady}{ad-bc}. \\
\end{array}
\end{equation} 
The local geometric meaning of the first integral 
$f=\mathfrak{c}_3 (Y)$ is 

\centerline{
$f(x,y)-f(x_0,y_0)= 
\left\{ \begin{array}{l}
\text{time of } Y \text{required to travel between the}\\
\text{trajectories of } X \text{ by } (x,y) \text{ and } (x_0,y_0)
\end{array}
\right\}$.}

\noindent 
The first integrals $f= \mathfrak{c}_3 (Y)$  are in general multivalued
(for example when $X$ has a source or sink). Therefore, the 
sheaf structure on $\mathcal{FI}(X)$ allow us 
 to use multivalued functions. 
\\
As an observation,
$\mathfrak{c}_2$ do not becoming the inverse operator to
$\mathfrak{c}_3$. 

\noindent 
The remarkable {\it Lie integrating factor $\mu$}, 
\cite{Ibragimov} p.\,267., determines the operator
\begin{equation}
\begin{array}{ccccccc}
\mathfrak{c}_5 : & Sym (X) & \longrightarrow  & \mathcal{IF}(X), \ \ \
& Y = c\del{}{x} + d\del{}{y} & \longmapsto & 
\mu = \frac{-1}{ad-bc}, \\
\end{array}
\end{equation}

\noindent
We propose the fourth operator as
\begin{equation}
\begin{array}{ccccccc}
\mathfrak{c}_4 : & \mathcal{IF}(X) & \longrightarrow  & Sym (X), \ \ \ 
& \mu & \longmapsto & 
Y= \left\{
\begin{array}{ll}
\frac{ 1}{ \Vert \nabla \mu \Vert ^2 } \nabla \mu  & \hbox{ if } X \hbox{ is Hamiltonian} \\
\frac{1}{\mu^2 div (X)} X_\mu  & \hbox{ if } X \hbox{ is non Hamiltonian.} \\
\end{array} \right.
\\
\end{array}
\end{equation}

\noindent 
If $X$ is a Hamiltonian vector field, then each integrating factor is
a first integral. 
However, $\mathfrak{c}_5$ and $\mathfrak{c}_4$ are not inverse one of the other.

\begin{example}
\begin{upshape}
In general, the domains where the first integrals, 
symmetries and 
integrating factors are $C^\infty$ do not coincide. 
The Lotka--Volterra vector field is

\centerline{$
X(x,y)=(\alpha x+\delta xy) \del{}{x} + (\beta y+\eta xy) \del{}{y},
\ \ \
\alpha, \beta, \delta, \eta \in \RR, \ 
\alpha, \beta \neq 0. 
$}

\noindent 
Let us consider 

\centerline{$
f(x,y) = \alpha \ln y - \beta \ln x + \delta y - \eta x 
\ \in \ \mathcal{FI}(X) 
\ \ \
\hbox{on } U_1 \doteq \RR^2 \backslash \{xy=0\}.
$}

\noindent
The operators $\mathfrak{c}_1$ and $\mathfrak{c}_2$
determine an integrating factor and a symmetry

\centerline{$
\begin{array}{c}
\mathfrak{c}_1 \big(f(x,y) \big) = \displaystyle -\frac{1}{xy}, 
\\
\mathfrak{c}_2 \big(f(x,y) \big) = \displaystyle
\frac{xy}{(\alpha + \delta y)^2 x^2 + (\beta + \eta x)^2 y^2}
\left (-(\beta y+\eta xy) \del{}{x} + (\alpha x+\delta xy) \del{}{y} \right)
\end{array}
$}

\noindent 
on domains $U_1 = \RR^2 \backslash \{xy=0\}$ 
and $U_2 \doteq \RR^2 \backslash 
\big( \overline{0} \cup (-\frac{\beta}{\eta}, 
-\frac{\alpha }{ \delta} )  \big)$. 
\end{upshape}
\end{example}

\section{Proof of Theorem 1.}

As motivation consider a naive question. Under what conditions the 
Hamiltonian and the gradient vector fields of a function commute? 
Let $J_0$ be the canonical complex structure on $\RR^2$. 

\begin{corollary} 
\label{conmutador-hamiltoniano-gradiente-armonico}
\begin{enumerate}
\item 
On $\CC \doteq (\RR^2, J_0)$
the following assertions are equivalent. 

\begin{enumerate}
\item 
The function
$V : \RR^2 \longrightarrow \RR$ is a harmonic.

\item
The Hamiltonian $X_V$ and the gradient vector field
$\nabla V$ of $V$ 
commute  on $\RR^2$ 
up to repara\-me\-tri\-zation 
by 

\centerline{
$\rho = (V_x ^2 + V_y^2)^{-1}$, 
\ \ \ 
thus $[\rho X_V, \rho \nabla V] = 0$.
}
\end{enumerate}

\item 
Moreover, for any Riemann surface 
$(\RR^2 \backslash \mathfrak{S}, J)$, 
the equivalence (i)--(ii) remains true
for $J$--harmonic functions,
where $\rho$ depends on $J$.  
\end{enumerate}
\end{corollary}

\begin{proof}
The function $V$ 
determines a pair of 1--forms and its dual vector fields, 

\centerline{$
\begin{array}{ll} 
\omega_X\doteq V_xdx+V_ydy & 
\omega_Y\doteq -V_ydx+V_xdy\\
X\doteq  \rho \nabla V
=\frac{1}{V_x^2+V_y^2}\left(V_x\del{}{x}+V_y\del{}{y} \right) 
\ \ \
&
Y \doteq  \rho X_V
=\frac{1}{V_x^2+V_y^2}\left(-V_y\del{}{x}+V_x\del{}{y} \right), 
\end{array}
$}

\noindent 
see L.\,V.\,Ahlfors \cite{Ahlfors1},  pp.\,162--163.  
Now, the equivalence (i)--(ii) is a routine computation. 
\end{proof}

\begin{proposition}
\label{coorrespondencia-basica}
\begin{enumerate}
\item 
On $\CC$ there exists a 
natural one to one correspondences between singular complex analytic
vector fields $\XX$, singular complex analytic 
1--forms $\omega$
and singular complex analytic maps $\Psi$ 
(probably multivalued but having single valued differential 
\footnote{These functions are called additively automorphic.}), as the diagram shows
\begin{equation}
\label{trilogia-complejo-analitica}
\begin{picture}(150,48)
\put(65,42){$\XX (z) = h(z)\del{}{z}$}
\put(-5,0){$\omega(z)=\frac{dz}{h(z)}$}
\put(137,0){$\Psi (z)=\int^z \omega$}

\put(60,3){\vector(1,0){75}}
\put(135,3){\vector(-1,0){82}}
\put(65,36){\vector(-3,-2){38}}
\put(27,11){\vector(3,2){38}}
\put(161,11){\vector(-3,2){38}}
\put(126,34){\vector(3,-2){38}}
\end{picture}
\end{equation}

\noindent  
the correspondence with $\Psi$ 
is up to additive constant.

\item
Moreover,
the following equalities hold 
\begin{equation}\label{dualidades}
\omega (\XX) = 1, 
\ \ \
d \Psi = \omega ,
\ \ \
\Psi_* \XX = \del{}{t} ,
\end{equation}
here $t= \tau+ \sqrt{-1} \sigma $ 
is the target variable of $\Psi$ and complex time of $\XX$.

\item 
For any Riemann surface
$(\RR^2 \backslash \mathfrak{P}, J)$,
the analogous correspondence
\eqref{trilogia-complejo-analitica}  
remains true.  
\hfill $\Box$ 
\end{enumerate}
\end{proposition}

The left arrow in \eqref{trilogia-complejo-analitica}
is implicit in Ahlfors \cite{Ahlfors1},  
pp.\,162--163, 
see also \cite{Mucino-Valero}, \cite{Mucino} and 
\cite{Alvarez-Mucino}.
The accurate application 
of \eqref{trilogia-complejo-analitica}
can be conducted as in the following possibilities
{\bf 1}--{\bf 4}, depending on the starting data.

\smallskip

\noindent
{\bf 1.}  {\it Let 
$X = u \del{}{x} + v\del{}{y} \in 
\mathfrak{X}^\infty (\RR^2 \backslash \mathfrak{S})$  
be a real vector field  
satisfying the Cauchy--Riemann equations.
}

\noindent 
Let
$
Y \doteq -v\del{}{x} + u\del{}{y}
$ 
be the rotated vector field, under
the canonical complex structure $J_0$, 
we obtain a complex analytic vector field
$$
\XX =\big( u+\sqrt{-1}v \big) \del{}{z} = \frac{1}{2} \big(X -\sqrt{-1} {J_0}X \big),
$$
see 
Kobayashi {\it et al.} \cite{Kobayashi-Nomizu} p.\,129,  , 
Prop. 2.11.  
By definition,
\begin{equation}\label{campos-parte-real-e-imaginaria}
X \doteq \Re{\XX} 
\ \ \ \hbox{and} \ \ \ Y \doteq \Im{\XX}
\end{equation} 
are {\it the real  and imaginary part of $\XX$}.
They commute 
and are linearly independent on  
$M_1 \doteq 
\RR^2 \backslash (\mathcal{Z}(X) \cup \mathfrak{S})$.

\noindent  
The dual frame of real 1--forms is

\centerline{$
\omega_X=\frac{1}{u^2 + v^2} (udx + vdy) \ , \ \
\omega_Y=\frac{1}{u^2 + v^2} (-vdx + udy),$
}

\noindent
satisfying 
\vspace{-3mm}
\begin{equation}\label{marco-real}
\begin{array}{ccc}
\omega_X(X)= 1 &  , & \omega_X(Y) =0, \ \\ 
\omega_Y(X) =0 & ,  & \omega_Y (Y)=1, 
\end{array}
\end{equation}
that is the first equation in (\ref{dualidades}).
The (multivalued) global flow box is
$$
\Psi(x,y)=
\big(U(x,y),V(x,y) \big)= 
\Big(\int^{(x,y)} \omega_X, \int^{(x,y)} \omega_Y \Big): 
\RR^2 \backslash \mathfrak{P} \longrightarrow \RR^2
$$
for $X$ and $Y$ respectively 
($U,\ V$ remain real analytic at the poles of $\XX$).

\noindent
The third equation in (\ref{dualidades}) assumes the real form
\begin{equation}\label{caja-de-flujo-real}
(U, V)_* X = \del{}{\tau} \ , \ \  
(U, V)_* Y= \del{}{\sigma}.
\end{equation}

\smallskip 

\noindent
{\bf 2.} {\it  Let
$\Psi=(U, V): \RR^2 \backslash \mathfrak{S} \longrightarrow \RR^2$ 
be a $C^\infty$  map, 
satisfying the Cauchy--Riemann equations.}

\noindent 
Considering a critical points 
$\Sigma(U,V) \doteq\{ (x,y) \in \RR^2 \backslash \mathfrak{S} 
\ \vert \ U^2_x +  U^2_y =0\}$. 
We get two canonically associated real vector fields
\begin{eqnarray}
\label{hamiltoniano-gradiente}
X &\doteq & 
\frac{1}{ V^2_x + V^2_y}\left(V_y \del{}{x} - V_x \del{}{y}\right) =
\frac{-X_{V}}{ V^2_x + V^2_y},\\
Y  &\doteq & 
\frac{1}{ U^2_x + U^2_y}\left(-U_y \del{}{x} + U_x \del{}{y}\right) =
\frac{X_U}{ U^2_x + U^2_y}
=
\frac{1}{V_x ^2 + V_y ^2} \left( V_x \del{}{x} + V_y\del{}{y}\right).
\nonumber
\end{eqnarray}   
The associated singular complex analytic vector field is

\centerline{$
\XX 
= \frac{1}{2} \Big(X + \sqrt{-1} {J_2} X \Big)= 
\Big(\frac{1}{U + \sqrt{-1}V}\Big)\del{}{z}.
$}

\noindent
Note that 
$J_2=J_0$, $Y$ is linearly independent with $X$ on 
$M_2 \doteq \RR^2 \backslash \big(\Sigma(U,V) \cup \mathfrak{S} \big)$ 
and $[X,Y]=0$.

\noindent
The dual frame of 1--forms is
$$
\omega_X = U_x dx + U_y dy =dU \ , \ \
\omega_Y = V_x dx + V_y dy =dV
$$
satisfying (\ref{marco-real}). Therefore (\ref{caja-de-flujo-real})
remains true in this case.

\smallskip

\noindent
{\bf 3.} {\it Let  
$(g, f): \RR^2 \backslash \mathfrak{S}  \to \RR ^2 $
be a $C^\infty$ map, 
which is a local diffeomorphism 
and a set $\mathfrak{S}$
with at most a finite number of accumulation points.
No Cauchy--Riemann conditions are required.}

\noindent 
The singular set as a map is 
$\Sigma (g,f) =\{ (x, y) \in \RR^2 \backslash \mathfrak{S} 
 \ \vert \ g_x f_y - f_x g_y =0 \}$.  
Using equation (\ref{dualidades}),
we get two canonically associated real vector fields
$$
X \doteq
(g,f)^*\del{}{t} =
\frac{1}{ g_x f_y - f_x g_y}\left(f_y \del{}{x} - f_x \del{}{y}\right) =
\frac{ -X_f}{ g_x f_y - f_x g_y},
$$
$$
Y \doteq
(g,f)^*\del{}{s}
= \frac{1}{ g_x f_y - f_x g_y}\left(-g_y \del{}{x} + g_x \del{}{y}\right) =
\frac{X_g}{ g_x f_y - f_x g_y}.
$$

\noindent 
Note that $Y$ is linearly independent with $X$ on 
$M_3\doteq \RR^2 \backslash (\Sigma (g,f) \cup \mathfrak{S})$
and clearly $[X, Y]=0$. 

\noindent
The dual frame of $1$--forms is 
$$
\omega_X = g_x dx + g_y dy =dg \ , \ \
\omega_Y = f_x dx + f_y dy =df.
$$

\noindent 
We regard to the canonical complex structure, say $J_3$, 
on the target $\{ (\tau, \sigma)\} 
=\{ \tau +  \sqrt{-1}\sigma \}$ 
of $(g,f)$ and the pull--back complex structure on the domain 
\begin{equation}\label{estructura-compleja}
J_3 (X) \doteq Y \ , \ \ J_3 (Y)\doteq -X .
\end{equation}

\noindent
Note that, a priori $J_3\neq J_0$, 
it is different from the canonical structure on $M_3$.
As a result, 
the map $(g,f): (M_3, J_3)\longrightarrow \CC $ 
becomes a local biholomorphism between 
Riemann surfaces, 
\cite{Kobayashi-Nomizu} p.\,115.   
Therefore, 
$$
\XX = \frac{1}{2} \big(X + \sqrt{-1}  J_3 X \big)
$$
\noindent is a singular complex analytic vector field
respect to $J_3$, see \cite{Kobayashi-Nomizu} p.\,122. 

By definition given a pair $(g,f)$ as above, 
$g$ is the {\it mate} of $f$.

\smallskip

\noindent
{\bf 4.}
{\it Let  
$X = a \del{}{x} + b\del{}{y} \in 
\mathfrak{X}^\infty (\RR^2 \backslash \mathfrak{S})$ 
be a  real vector field 
admitting a second one
$Y = c \del{}{x} + {\tt d} \del{}{y}$ that commutes, 
this is $[X,Y] \equiv 0$. 
}

\noindent 
The $C^\infty$ singular set is 
$Sing(X, Y) =\{ (x,y) \in \RR^2
 \ \vert \ a{\tt d}-bc=0 \}$.  
We consider 

\centerline{$M_4\doteq \RR^2 \backslash (Sing(X,Y)\cup \mathfrak{S})$}

\noindent 
and define an adapted complex structure $J_4$ as 
\begin{equation}\label{estructura-compleja-2}
J_4 X \doteq Y \ , \ \ J_4 (Y)\doteq -X .
\end{equation} 
Then the pair $(M_4, J_4)$ is a Riemann surface.  
The dual frame of 1--forms  is

\centerline{$
\omega_X=\frac{1}{a{\tt d}-bc} ({\tt d}dx - cdy) \ , \ \
\omega_Y=\frac{1}{a{\tt d}-bc} (-bx + ady), 
$}

\noindent 
satisfying (\ref{marco-real}). 
They are closed 1--forms by the integrability hypothesis. 
There are two (probably multivalued) first integrals 
\begin{equation}\label{formula-integral-multivaluada}
\big(g(x,y), f(x,y)\big)= 
\Big( \int^{(x,y)} \omega_X , \int^{(x,y)} \omega_Y \Big): 
\RR^2 \backslash \big(Sing(X,Y) \cup \mathfrak{S}\big) \longrightarrow  \RR^2 
\end{equation}

\noindent 
of $Y$ and $X$ respectively, such that
the map $(g, f): (M_4, J_4) \to \CC $ is a local biholomorphism,
see \cite{Kobayashi-Nomizu} p.\,122. 
Therefore,
$$
\XX = \frac{1}{2} (X + \sqrt{-1}  {J_4} X)
$$
\noindent 
is a complex analytic vector field,
respect to ${J_4}$.  
The real form of (\ref{dualidades}) remains true.

\noindent  
$X$ is the infinitesimal generator of a locally free 
$(\RR^2, +)$--action on $M_4$.

\smallskip

We summarize the diagram \label{trilogia-compleja} 
and the possibilities {\bf 1}--{\bf 4} as follow. 

\begin{proposition}\label{integrabilidad-caja-flujo-global}
\begin{enumerate}
\item 
On the respective non
singular loci $M_\iota$, $\iota \in 1, \ldots , 4$,
there exists a natural one to one correspondence 
\begin{equation}\label{correspondencia-real}
\begin{picture}(350,58)
\put(0,-4){$\Psi (z) = \int^z \frac{d\zeta}{h(\zeta )} \ on \ (M_2, J_0)$}
\put(160,-4){$\big( g(x,y),f (x,y) \big) \ on  \ (M_3, J_3),$}
\put(0,50){$\XX (z) = h(z)\del{}{z} \ on  \ (M_1, J_0)$}
\put(160,50){$X(x,y), \ Y(x,y), \ [X,Y] \equiv 0 \ on  \ (M_4, J_4)$}

\put(134,-1){\vector(1,0){23}}
\put(154,-1){\vector(-1,0){23}}

\put(134,53){\vector(1,0){23}}
\put(154,53){\vector(-1,0){23}}

\put(58,40){\vector(0,-1){30}}
\put(58,11){\vector(0,1){31}}

\put(208,40){\vector(0,-1){30}}
\put(208,11){\vector(0,1){31}}

\end{picture}
\end{equation}

\vskip0.2cm

\noindent 
here in the left column, $z$ is a complex variable respect to the 
suitable adapted complex structure.

\item 
The complex analytic $\XX$,  
real $\Re{\XX}$ and 
Hamiltonian $X_f$ vector fields 
in (\ref{correspondencia-real}) 
admit  

\begin{enumerate}
\item 
global flow box $\Psi=g + \sqrt{-1}f$,

\centerline{$\Psi_* \XX = \del{}{t}, 
\ \ \
(g, f)_* X_f = \del{}{\tau}$,} 
 
\item
and the adapted complex structures $J_\iota$ 
are such that 

\centerline{$ X = \Re{\XX} = X_f$.}
\end{enumerate}
\end{enumerate}
\end{proposition}

\begin{proof}
By simple
inspection, we start with the respective
non singular object on  $(M_\iota, J_\iota)$ 
and calculate the other three objects.   
\end{proof}

Let us make some observations about (\ref{correspondencia-real}).
Two vector fields $X$, $Y \in\mathfrak{X}^\infty (M_4)$ 
are orthogonal 
and $\vert X \vert = \vert Y \vert $ if and only if 
the corresponding  $J_4$ 
is the canonical complex structure $J_0$.

\begin{remark}[Non uniqueness of 
the global flow box map $(g,f)$ in 
\eqref{ecuacion-de-caja-de-flujo-global},  
\eqref{correspondencia-real}]
\label{la-caja-de-flujo-no-es-unica}
\begin{upshape}
The group of $C^\infty$ diffeomorphisms 
of $\RR^2 = \{ (\tau, \sigma )\}$
preserving the vector field $\partial / \partial \tau$ 
is generated by diffeomorphisms of type;
$\phi_\iota  (\tau, \sigma)= (\tau + h_\iota( \sigma ),  \sigma )$ 
\,  called 
{\it shear} or {\it Jonqui\`ere} 
maps and 
$\phi_j(\tau, \sigma)= (\tau , h_j(\sigma))$,  
here $h(\sigma) \in \hbox{ Diff}^\infty(\RR, \RR)$. 

\noindent 
Hence, the global flow box map 
$(g,f)(x, y)=(\tau, \sigma )$ of $X$ is far from being unique, 
{\it i.e.}
the transversal structure of $X$ is non canonical. 
\end{upshape}
\end{remark}

\smallskip

\noindent
{\it Proof of Theorem 
\ref{teorema-1-introduccion-diccionario}} follows by simple inspection of
Proposition \ref{integrabilidad-caja-flujo-global}.

\section{Canonical regions for complex vector fields}

\begin{figure}[h!]
\begin{center}
\scalebox{0.48}{\includegraphics{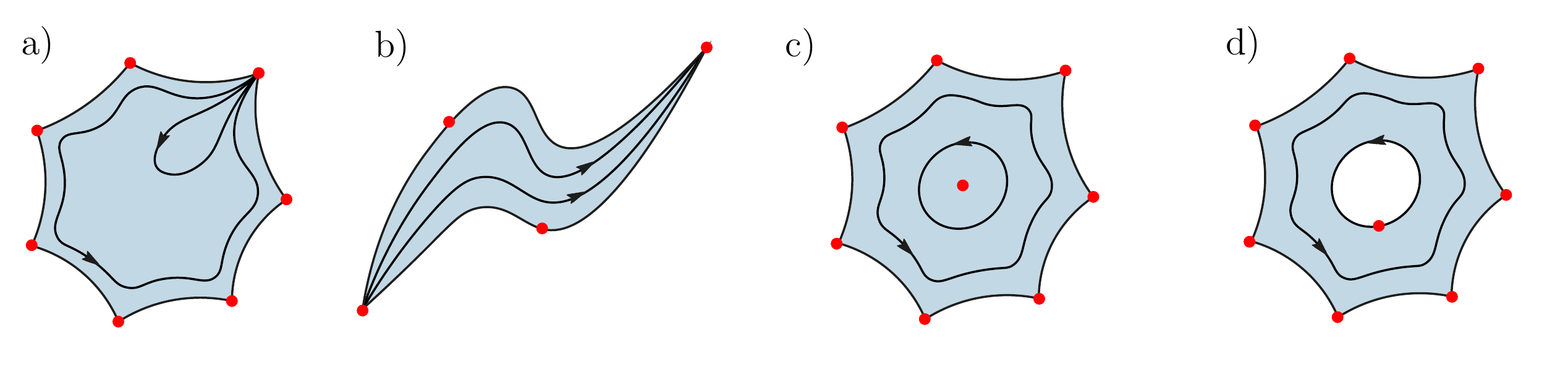}}
\caption{
Real phase portraits of 
canonical regions of singular complex analytic vector fields $\XX$,
in their boundaries, pieces of trajectories
and singular points (in red) appear. }
\label{regiones-Markus-Benzinger}
\end{center}
\end{figure}

\begin{proposition}
There exists a flat Riemannian metric $g_\XX$ associated
to $\XX$ on $(\RR \backslash \mathfrak{P}, J)$, 
such that the real trajectories of $\Re{\XX}$ are unitary geodesics. 
\hfill $\Box$ 
\end{proposition}
\begin{proof}
For the proof see \cite{Mucino-Valero}, \cite{Mucino}, or/and 
\cite{Jenkins},
\cite{Strebel} for the quadratic differentials point of view. 
\end{proof}
\begin{definition}\label{regiones-de-Fatou} 
\begin{upshape}
\begin{enumerate}
\item 
The {\it open canonical regions} of  are pairs 
(domain \& holomorphic vector field)  
as follows 
\begin{equation}\label{ecuaciones-regiones-Fatou}
\begin{array}{rcr}
\hbox{\it half plane }
\mathscr{H}=
\big( \mathbb{H}^2,\frac{\partial}{\partial z}\big),
&&
\hbox{\it  strip }  
\mathscr{S}=
\big( \{ 0 < \Im{z} < h \},
\frac{\partial}{\partial z}\big),
\\
&& \vspace{-.2cm}
\\
\hbox{\it half cylinder }  
\mathscr{C}=
\Big(\Delta_1, 
\dfrac{-2 \pi i z}{r}\frac{\partial}{\partial z} \Big),
&&
\hbox{\it annulus }  
\mathscr{A}=
\Big(\Delta_R \backslash 
\overline{\Delta}_1, \dfrac{-2 \pi iz}{r}
\frac{\partial}{\partial z}\Big),
\\
&& \vspace{-.3cm}
\end{array}
\end{equation}
\noindent  
here  $\mathbb{H}^2 = \{ \Im{z} > 0 \}$ is the open half plane 
and $\Delta_R = \{ \vert z \vert < R \}$
is an open disk.
 
\item
Given $\XX$ on $(\RR^2 \backslash \mathfrak{P}, J)$, 
a pair $(\mathcal{U}, \XX)$ is a {\it canonical region}
of $\XX$ 
when it is holomorphically equivalent to one element in  
\eqref{ecuaciones-regiones-Fatou} and it is maximal.  
See Figure \eqref{regiones-Markus-Benzinger}.
\end{enumerate}
\end{upshape}
\end{definition}

The canonical regions are 
$\Re{\XX}$--invariant, in particular their  real
trajectories $\{ \tau \mapsto z(\tau) \} $ are 
well defined for all real time each canonical region.
The boundaries of the canonical regions are geodesics in the 
metric $g_\XX$. 
The factor $-2\pi i z/ r$ in \eqref{ecuaciones-regiones-Fatou} 
makes the geodesic boundaries of $g_\XX$--lenght $r> 0$, in the case of half cylinder
and annulus.

\begin{definition}\label{separatriz-holomorfas-definicion}
\upshape
Let $\XX$ be a singular complex analytic vector field
on $(\RR \backslash \mathfrak{P}), J)$, 
the {\it separatrix skeleton 
$\Gamma(\XX)= \{ z_\vartheta (\tau) \}$
of $\XX$} is the union of its 
$\RR$--incomplete trajectories. 
\end{definition}

\begin{example} \label{ejemplo-polo}
{\it A pole of order $k$.}
\begin{upshape} 
Let $\XX$ be singular 
complex vector field on $\CW$, having a pole of order $k \geq 1$, 
the diagram from Equation \eqref{correspondencia-real} is
\begin{center}
\begin{picture}(300,67)
\put(11,-4){$\Psi(z) = \int^z \zeta^k d \zeta $}
\put(70,30){\vector(0,1){18}}
\put(70,30){\vector(0,-1){18}}
\put(27,57){$\XX(z) = \frac{1}{z^k} \del{}{z}$}
\put(110,59){\vector(1,0){20}} 
\put(110,59){\vector(-1,0){20}} 
\put(128,56){$
\begin{array}{c}
X(x, y) = \Re{\frac{1}{z^k}} \del{}{x} + \Im{\frac{1}{z^k}}\del{}{y} \\
\ \ Y (x, y) = - \Im{\frac{1}{z^k}}\del{}{x} + \Re{\frac{1}{z^k}} \del{}{y} 
\end{array}$}

\put(110,-2){\vector(1,0){20}} 
\put(110,-2){\vector(-1,0){20}}

\put(200,25){\vector(0,1){14}}
\put(200,25){\vector(0,-1){14}} 
\put(135,-4){$(g, f)= \Big( \Re{\frac{z^{k+1}}{k+1}}, 
\Im{\frac{z^{k+1}}{k+1}} \Big).$}
\end{picture}
\end{center}

\noindent 
The canonical decomposition of $\XX$ is

\centerline{$
\CW \backslash \Gamma (\XX)= 
\bigcup_{\alpha=1}^{2k+2}
\big( {\HH}^2, \del{}{z}\big)_\alpha,
$}

\noindent
having complete separatrix skeleton 

\centerline{$ \Gamma (\XX) = \left\{ \Im{z^k} = 0 \right\} \cup \left\{ \infty \right\}
\subset \CW. 
$}

\noindent 
Looking at $\infty, 0 \in \CW$, 
the germs\footnote{The change
of coordinates for $\CW$ is 
$z \mapsto w=1/z$.}
of $\XX$ and admissible words are

\centerline{$
\big( 
(\CC_z, 0 ),
\frac{1}{z^k}\del{}{w} \big)
\longleftrightarrow
\underbrace{H \cdots H}_{2k +  2},
\ \ \
\big( 
(\CC_w,0), w^{k-2}\del{}{w} 
\big)
\longleftrightarrow
\underbrace{E \cdots E}_{2(k-2)},
$ for $k\geq 3.
$}

\noindent 
Here $H$, $E$ means the hyperbolic, elliptic angular sectors of $\Re{\XX}$, 
see Table 1 and Figure \ref{nuevo-album-2}.a--b. 
\end{upshape}
\end{example}

\begin{figure}[h!]
\begin{center}
\scalebox{0.40}{\includegraphics{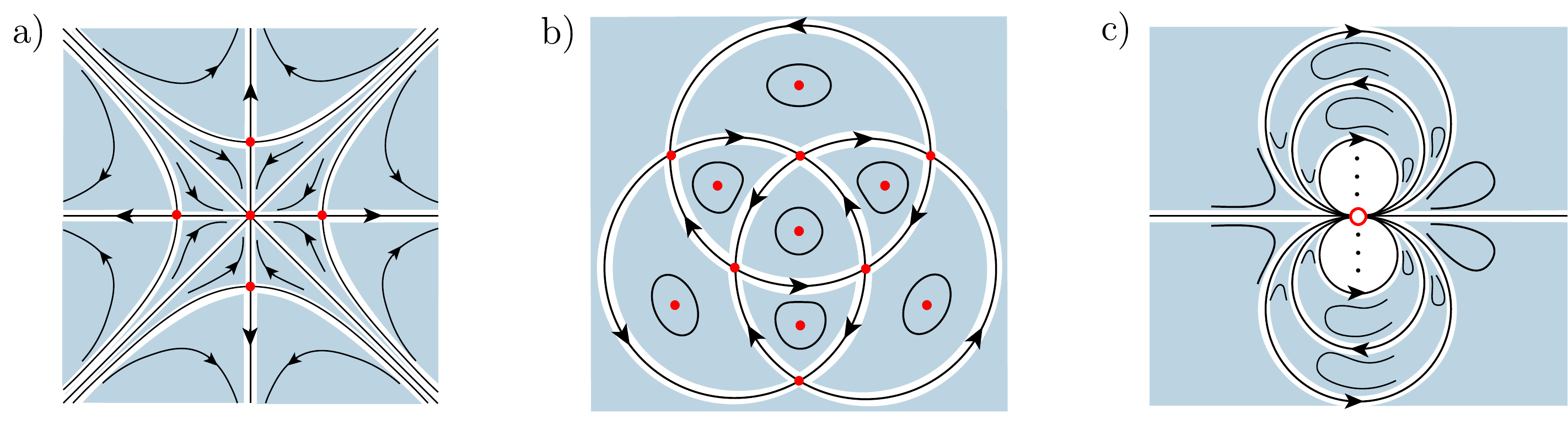}}
\caption{
Decomposition in canonical regions without accumulation of singular points (in red) 
of $\XX$ at $\infty \in \CW_z$. (a) and (b) are rational vector fields, (c) is the
complex exponential vector field (the small red circle denotes $\infty$).}
\label{pegados1}
\end{center}
\end{figure}

\begin{example} \label{ejemplo-solecito-real}
\begin{upshape}
Consider the complex rational vector field as follows 

\centerline{
$\XX(z) =\frac{z}{z^4-1} \del{}{z}$, 
\ \ \
$X(x,y) \doteq \rho(x,y) \Re{\frac{z}{z^4-1} \del{}{z}  }$, 
}

\noindent
for suitable $\rho$, 
the last
is $C^\infty$ on $\RR^2$. 
Using  Figure \ref{pegados1}.a, we note that
$\CC \backslash \Gamma(\XX) $
is a union de eight strips and eight half planes.  
\end{upshape} 
\end{example}

\begin{example} \label{ejemplo-exponencial}
{\it The complex exponential vector field on $\CW$. }
\begin{upshape} 
The diagram is
\begin{center}
\begin{picture}(320,75)
\put(7,0){$\Psi(z) = \int^z \e^{-\zeta} d \zeta $}
\put(65,34){\vector(0,1){20}}
\put(65,34){\vector(0,-1){20}}
\put(27,61){$\XX(z) = \e^z \frac{\partial}{\partial z}$}
\put(112,63){\vector(1,0){20}} 
\put(112,63){\vector(-1,0){22}} 
\put(130,60){$
\begin{array}{c}
X(x, y) = \e^x cos (y) \frac{\partial}{\partial x} + \e^x sin(y)\del{}{y} \\
\ \ Y (x, y) = - \e^x sin (y) \del{}{x} + \e^x cos (y) \del{}{y}
\end{array}$}

\put(112,2){\vector(1,0){20}} 
\put(112,2){\vector(-1,0){22}}

\put(182,29){\vector(0,1){16}}
\put(182,29){\vector(0,-1){16}} 
\put(137,0){$(g, f)= (\Re{\Psi}, \Im{\Psi})$.}
\end{picture}
\end{center}

\noindent 
The vector field 
$X=\Re{\XX}$ have an infinite number of horizontal Reeb's components
on $\CC$, Figure \ref{pegados1}.c.
The decomposition in canonical regions is

\centerline{$
\CW \backslash \Gamma (\XX) = 
\bigcup_{\alpha=-\infty}^\infty
\big( \HH^2, \del{}{z}\big)_\alpha,
$}

\noindent
having separatrix skeleton 

\centerline{$
\Gamma (\XX) = \{ \Im{z} = k \pi \ \vert \  k \in \ZZ\} \cup \{ \infty \}. $}

\noindent 
Looking ${\XX}$ at the point $\infty \in \CW$, 
the germ and admissible word are

\centerline{$
\big( 
(\CC_w, 0),
-w^2\e^{-w}\del{}{w} 
\big)
\longleftrightarrow
\ent{}{} \ent{}{}.
$}

\noindent 
\noindent 
Here $\ent{}{}$ means the entire sector of $\Re{\XX}$, 
see Equation \eqref{sectores-holomorfos} and Figure \ref{nuevo-album-2}.d-e, 
the entire sectors
come from $\XX$ on 
$\{z \ | \ \Im{z} > 0\}$ and 
$\{z \ | \ \Im{z} < 0\}$.
\end{upshape}
\end{example}

%
%

%
 
%

%
%

\begin{example}
{\it An infinite number of isochronous centers.}
\begin{upshape}
Recalling that, 
a simple zero with pure imaginary linear part
determines an isochronous center,  we consider
\begin{center}
\begin{picture}(330,75)
\put(0,-4){$\Psi(z) = \int^z \frac{1}{i \, sin(\zeta)} d \zeta $}
\put(65,30){\vector(0,1){20}}
\put(65,30){\vector(0,-1){20}}
\put(9,57){$\XX(z) = i\, sin(z) \del{}{z}$}
\put(115,59){\vector(1,0){20}} 
\put(115,59){\vector(-1,0){20}} 

\put(131,56){$
\begin{array}{c}
X(x, y) = \Re{ i \, sin(z)} \frac{\partial}{\partial x} + 
\Im{ i \, sin(z)} \del{}{y} \\
\ \ Y (x, y) = - \Im{ i \, sin(z)} \del{}{x} + \Re{ i \, sin(z)} \del{}{y}
\end{array}$}

\put(115,-2){\vector(1,0){20}} 
\put(115,-2){\vector(-1,0){20}}

\put(183,25){\vector(0,1){16}}
\put(183,25){\vector(0,-1){16}} 
\put(139,-4){$\Big( \Re{ i \, csc(z) \, cot(z)}, \Im{ i \, csc(z) \, cot(z)} \Big),$}
\end{picture}
\end{center}

\noindent 
see Figure \ref{sectores-esenciales-nuevos}.b. 
The canonical decomposition is 

\centerline{$
\CW \backslash \Gamma (\XX)= 
\bigcup_{\varsigma=1}^\infty
\big
(\Delta_1, \pm iz \frac{\partial}{\partial z} 
\big)_\varsigma ,$}

\noindent
having an infinite number of half cylinders (isochronous centers).
\end{upshape}
\end{example}

\section{Proof of Theorem 2}
\label{prueba-de-teorema-2}

\smallskip 
\noindent
{\it Proof assertion $(1)$ in Theorem
\ref{teorema-2-introduccion-Fatou-decomposition}. }
The main technical result for the proof is the following. 

\begin{lemma}[\cite{Strebel} 
p.\,57, 
\cite{Mucino-Valero},  
\cite{Mucino}; 
Isometric glueing]
\label{pegado-isometrico}
Let $(\overline{N}_1, \XX_1)$, $(\overline{N}_2, \XX_2)$ be two canonical regions
and let 
$z_1 (\tau)\subset \partial \overline{N}_1$, $z_2 (\tau) \subset \partial \overline{N}_2$ 
be segments 
in trajectories of $\Re{\XX_1}$ and $\Re{\XX_2}$, 
having the same length.
The isometric glueing of them 
along these geodesic boundary preserving the orientation of 
$\Re{\XX_1}$ and $\Re{\XX_2}$,  
is well defined, and provides a new flat surface (a Riemann surface structure)
on $\overline{N}_1 \cup \overline{N}_2$ arising from a new 
complex analytic vector field
$\XX$.
\hfill $\Box$
\end{lemma}

By hypothesis $M$, is
obtained by the paste of closed canonical regions 
$$
\left( 
\overline{\mathbb{H}}^2, \del{}{z} \right),\
\Big( \{ 0 \leq \Im{z} \leq h \}, \frac{\partial}{\partial z}\Big),\
\left( \overline{\Delta}_1, \frac{-2\pi iz}{r} \del{}{z}\right), \
\left( \overline{\Delta}_R \backslash \Delta_1, 
\frac{-2\pi iz}{r} \del{}{z} \right), 
$$

\noindent 
$h,\, r \in \RR^*$. 
Here the paste uses isometries that preserve the orientation 
of the real trajectories in their boundaries.
The assertion follows from the Corollary \ref{pegado-isometrico}.

\smallskip 
\noindent
{\it Proof assertion (2) in Theorem
\ref{teorema-2-introduccion-Fatou-decomposition}. }
By hypothesis, the incomplete trajectories of $\Re{\XX}$
are locally finite in $(M, J)$. 
Then, we remove the separatrix skeleton, thus  
$M \backslash \Gamma(\XX)$ is an finite or infinite union of open sets
invariant under the flow of $\Re{\XX}$. 
It is well known that interior of the open sets are 
necessarily as in Definition \ref{regiones-de-Fatou}. 
Thus, the Riemann surface $(M, J)$ 
has a decomposition 
\begin{equation} \begin{array}{r}
M=
\bigcup_{{\tt a}} \left( 
\overline{\mathbb{H}}^2, \del{}{z} \right)
\bigcup_{{\tt b}}
\left( \{ 0 \leq \Im{z} \leq h_{\tt b} \}, \del{}{z} \right)
\bigcup_{{\tt c}}
\left( \overline{\Delta}_1, \frac{2\pi iz}{r_{\tt c}} \del{}{z}\right)
\bigcup_{\tt d}
\left( \overline{\Delta}_{R_{\tt d} } \backslash \Delta_1, 
\frac{ 2\pi i z }{r_{\tt d} } \del{}{z} \right).
\end{array}
\end{equation}
 
\noindent 
The Theorem \ref{teorema-2-introduccion-Fatou-decomposition} is done.

\begin{example}{\it 
Singular complex analytic vector fields
having an infinite decomposition in canonical regions.}
\begin{upshape}
In Figure \ref{sectores-esenciales-nuevos}
drawing (a) is the complex exponential, (b) is  $i sin(z) \del{}{z}$.    
The (c)--(f) use Theorem 2; they illustrate accumulation of poles and
zeros to an essential singularity at $\infty$. 
\begin{figure}[h!]
\begin{center}
\scalebox{0.40}{\includegraphics{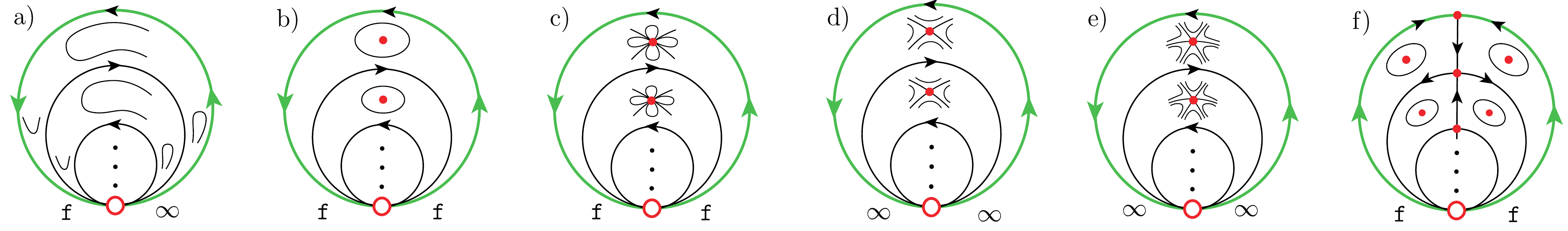}}
\caption{
Singular complex analytic vector fields $\XX$ on a half plane, 
with an infinite number of canonical pieces, 
accumulation of singular points (in red)
and infinite number of real incomplete trajectories of $\Re{\XX}$
at $ \infty \in \CW_z$.}
\label{sectores-esenciales-nuevos}
\end{center}
\end{figure}

\end{upshape}
\end{example}

Some results  on the set of incomplete trajectories 
$\Gamma(\XX)$ are in order. 

The local analytic normal forms for poles and zeros of
meromorphic complex analytic vector fields, 
Table \ref{tabla-polos-zeros}, is well known, 
here in the right column
$H,\, E, \, P$ means  topological hyperbolic, elliptic and parabolic 
topological sectors for $\Re{X}$.

\begin{table}
\caption{Local analytic normal forms for poles and zeros of $\XX$.}\label{tabla-polos-zeros}
\begin{center}
\begin{tabular}{|c|l|c|c|c|}
\hline
normal form & \ analytic  invariants & 
flat metrics& $\Psi(z)$ & top. invariants
\\
of  $\XX$ & \ \ \ \ order \& residue &
$g_\XX$ & & words 
\\
\hline
\hline
$\frac{\partial}{\partial z}$ & $\hspace{.8cm} 0 \hspace{.8cm}0$ & 
cone angle $2\pi$
& $z$  & $\underbrace{H \cdots H}_{2}$
\\
\hline
 $\frac{1}{z^{k}}\frac{\partial}{\partial z}$  & $ -k \leq -1 \hspace{.6cm} 0 $  
& cone angle $(2k + 2)\pi$ 
& $\frac{z^{k+1}}{k+1}$  & $\underbrace{H \cdots H}_{2k+2}$
\\
\hline
  $\lambda z \del{}{z}$ & \ \ \ \ $s=1 \ \ \ \ \ \lambda \in i\RR^* $ & 
$\mathbb{S}^1_{1\pi \vert \lambda \vert } \times (0, \infty)$ & $\lambda log(z)$ & $C$ \\
\vspace{-.3cm}
&&&& \\
\hline
 $\lambda z \del{}{z}$ & \ \ \ \ $s=1 \ \ \ \ \  \lambda \in \CC \backslash \RR$ &
$\mathbb{S}^1_{1\pi \vert \lambda \vert } \times (0, \infty)$
& $\lambda log(z)$ & $ P $ \\
\vspace{-.3cm}
&&&& \\
\hline
\vspace{-.5cm}
&&&& \\
 $\frac{z^{s}}{1 - \lambda z^{s-1}}\del{}{z}$  & \ \ \ \ $2 \leq s  \ \ \ \ \ \lambda \in \RR$ &
$(2s-2) \ \hbox{copies} \ ( \overline{\HH}^2, \infty)$ 
& $\frac{1}{(1-s)z^{s-1}} + \lambda log (z)$ & 
$\underbrace{E \cdots E}_{2s-2}$      \\
\hline
\vspace{-.5cm}
&&&& \\
&&$(2s-2) \ \hbox{copies} \ ( \overline{\HH}^2, \infty)$ & & \\ 
$\frac{z^{s}}{1 - \lambda z^{s-1}}\del{}{z}$ & \ \ \ \  $2 \leq s \ \ \ \ \ \lambda \in \CC \backslash \RR$ & 
and a strip
& $\frac{1}{(1-s)z^{s-1}} + \lambda log (z)$ & $\underbrace{E \cdots E}_{2s-2}P $  \\
\hline
\end{tabular}
\end{center}
\end{table}

\begin{corollary}
\label{caracterizacion-racional}
Let $\XX$ a singular complex analytic vector field on a 
simply connected  Riemann surface 
$(M, J)$ as above. 
The following assertions are equivalent. 

\noindent 
1) 
$(M, J) =\CC$ or $\CW$ and in any case $\XX$ is (the restriction of) a 
rational vector field on $\CW$. 

\noindent 
2) For each rotated vector field 
$\e^{i \theta} \XX$, its set of incomplete real trajectories $\Gamma(\e^{i \theta}\XX)$ is  
finite.

\noindent 
3) The decomposition of $\XX$ on $(M, J)$ has a finite number of canonical 
regions and no finite trajectory gap, {\it i.e.} a segment of geodesic in the boundary 
of a canonical region that is not identified in $(M,J)$.  
\end{corollary}

\begin{proof}
Assume Table.1 and let $\XX$ be a complex rational vector field on $\CW_z$. 
The incomplete trajectories are the separatrices at poles, hence they are finite 
number 
for each rotated vector field $\Re{\e^{i \theta} \XX}$, we perform the assertion (2). 
If we assume (2), then the separatrices are a finite number, and the decomposition 
in canonical pieces is finite. Moreover, the poles and zeros are conformal 
punctures of the complex structure $(M, J)$, this is no gap trajectory can appear; 
see Corollary 1 and its proof in \S\,\ref{seccion-tipo-conforme}. 
We leave the converse assertions for the reader.   
\end{proof}

\begin{figure}[h!]
\begin{center}
\scalebox{0.3}{\includegraphics{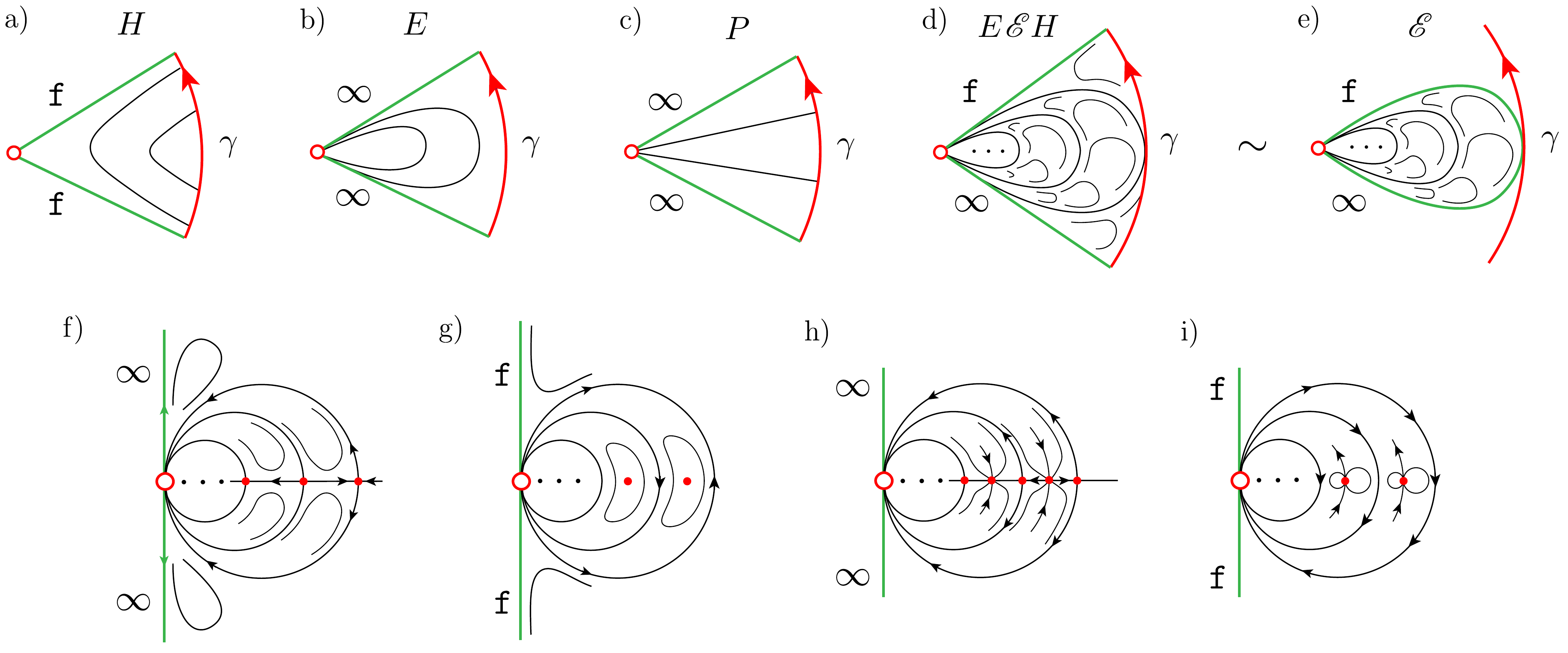}}
\caption{
Hyperbolic, elliptic and parabolic sectors are sketched in (a)--(c). 
Entire sector in (d)--(e). Other interesting sectors in (f)--(i).
The signs ${\tt f}$ and $\infty$ in the boundaries,
describe the {\it finite} or {\it infinite}
time of the boundary separatrix  to reach at the singularity at the vertex.}
\label{nuevo-album-2}
\end{center}
\end{figure}

Moreover, a local version of Theorem 2 provides
new non isolated singularities of complex analytic  vector fields, we give
a very brief description. 
The topological
{\it hyperbolic, elliptic and parabolic sectors for real vector fields}
of vector fields germs are illustrated in 
Figure \ref{nuevo-album-2}.a--c. 
Moreover they are analytic germs and suitable flat metrics
as follows.

The {\it
germs of singular complex analytic vector fields 
$\big((\CC, q), \mathbb{X} \big)$ on angular sectors} 
are as follows:
 \vspace{-.3cm}
\begin{equation}
\label{sectores-holomorfos}
\begin{array}{lll}
\hbox{a  \emph{isochronous center at } } q=0
& &
C=\left(  (\CC, 0), \ \frac{iz}{r} \del{}{z} \right), \
r \in \RR^+, 
\\
\hbox{a  \emph{hyperbolic sector at } }q=0
& &
H=\left( \big( \overline{\HH}^2, 0 \big),\ \del{}{z}\right),
\\
\hbox{an  \emph{elliptic sector at } } q=\infty
& &
E=\left( \big( \overline{\HH}^2, \ \infty \big), \del{}{z} \right),
\\
\hbox{a  \emph{parabolic sector, at } }
q=+\infty 
& &
P= \left( \left( \{  0 \leq \Im{z}  \leq h\}, 
\pm \infty \right), \ \del{}{z}\right), \ h \in \RR^+,
\\ 
\hbox{a \emph{1--class entire sector at }} q=\infty 
& 
\ \ \  
&
\ent{}{1}=
\left( \big( \overline{\HH}^2, \infty \big), \  
\e^{z} \del{}{z} \right). 
\end{array}
\end{equation}

In order to recombine the angular sectors in Figure \ref{nuevo-album-2}
for perform new singularities,
we consider as combinatorial information a  \emph{finite cyclic word}
\begin{equation}\label{palabra}
\mathcal{W}=
W_1 \cdots  W_k \ , \ \ \ W_\iota \ \hbox{ in the Figure 4 alphabet}.   
\end{equation}
The following conditions are satisfied by $\mathcal{W}$. 

\noindent 
Each $W_\iota$ can be interpreted also as a Riemannian manifold and 
an angular vector field germ $((\mathcal{A}_\iota,q_\iota ), \XX_\iota) \cong  W_\iota$.

\noindent 
The geodesic boundary of $W_\iota$ 
has orientation and time to reach the singular point 
denoted ${\tt f}$ when is finite or $\infty$, in our figures. 
We consider cyclic words in the sense that

\centerline{ $W_{k+1} \doteq W_1$.}

\noindent
In order to perform the geometric paste 
of the angular sectors, we require the 
following additional rule;
if the orientation and 
the time ${\tt t}$ or $\infty$ coincide
\\
then the two trajectories can be pasted together.
We have a new version of the result in \cite{Alvarez-Mucino} p.\,167:

\begin{corollary}
\label{de-palabras-a-germenes-de-campos-complejos}
\begin{enumerate}
\item 
The paste as above of a finite number of angular sectors 
$\mathcal{W}$
determines a germ of singular complex analytic vector field 
$\big( (\CC,0), \XX \big)$ with singularity at $0$.

\item
The singularity $0$ is a pole or zero of $\XX$ 
if and only if
a center,
a finite number of hyperbolic, elliptic and/or parabolic 
sectors appear at $0$, recall  Table 1.
\hfill 
$\Box$
\end{enumerate}
\end{corollary}

In the above, the vertex of the angular sectors is a
conformal puncture, hence the singularity is a zero, a pole or an essential
singularity of $\XX$.
The alphabet in Figure 4 is far from begin complete, however it shows the
wealth of the theory.

\section{On the conformal type problem}
\label{seccion-tipo-conforme}

\begin{figure}[h!]
\begin{center}
\scalebox{0.4}{\includegraphics{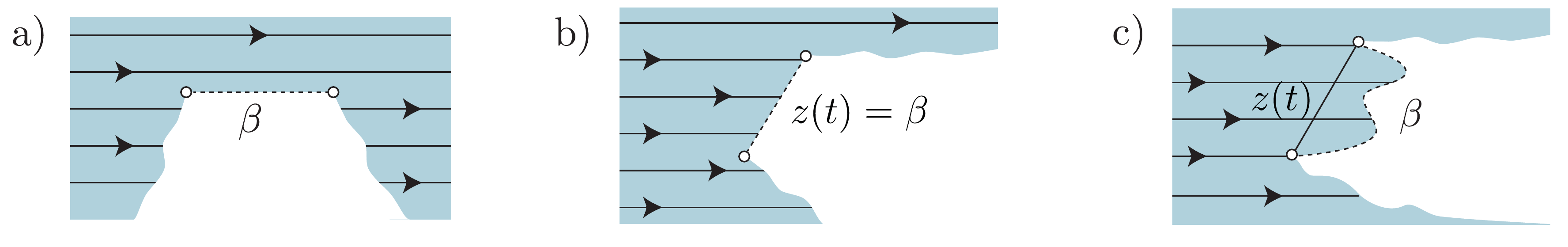}}
\caption{
A finite trajectory gap means that the ideal boundary of the metric $g_\XX$ 
(under a local flow box) is $\beta$. 
} 
\label{3-fig}
\end{center}
\end{figure}

\begin{definition}\label{gap}
\begin{upshape}
Let $J$ be complex structure on $M \subset \RR^2$ such that $X = \Re{\XX}$ is the
real part of a complex analytic vector field $\XX$ on $(M,J)$.
$\XX$ has a {\it finite trajectory gap} at the ideal boundary $\infty$
if there exists a local (holomorphic) flow box
$$
\Psi:(U\subset M,J)\longrightarrow \CC, \ \text{such that}\ 
\lim_{U\ni(x,y) \to \infty} \psi(x,y)=\beta\subset \CC,
$$
where 
$(x,y) \to \infty$ means that $(x,y)$  tends to the ideal boundary of $M$, and  
$\beta$ is a simple path in $\CC$ different from a point, see
Figure \ref{3-fig}.
\end{upshape}
\end{definition}

{\it Proof of Corollary \ref{corolario-2-introduccion-tipo-conforme}.}
Using the vector fields, we provide elementary arguments, 
compare with \cite{Ahlfors2} Ch.\,10.
Let $X$ be a vector field on a simply connected $(M,J)$ having a finite trajectory gap.

\noindent
Case 1. The path $\beta\in \CC$, gap is a segment of trajectory of $\XX$.
{\it i.e.} 
there exists in $\big((\RR^2, J), X \big)$ a holomorphic local flow box
$\Psi: U \subset \RR^2 \to  (0, \epsilon) \times (0, i \epsilon) \subset \CC$ such that 
$\Psi^{-1} (x + i0)$ is in the ideal boundary of $\RR^2$, 
Figure \ref{3-fig}.a. 

\noindent 
We proceed by contradiction, let $\pi : \CC \to (\RR ^2, J)$ the uniformization 
of the adapted complex structure to $X$. 
The map $\pi^{-1}$ 
sends the ideal boundary of $\RR^2$ to $\infty \in \CW$. 

\noindent 
Then the composition 
$(\frac{1}{z}) \circ \pi \circ \psi :(0, \epsilon) \times (0,i \epsilon)  \to  \CC$
is a biholomorphic map, continuous and real valued in $\{ x +i 0 \}$. 
Using the reflection principle, we can 
extend $(\frac{1}{z}) \circ \pi \circ \psi$ to a holomorphic
$\phi :(0, \epsilon) \times (-i\epsilon ,i \epsilon)  \to  \CC$, 
defining $\phi(x+iy)\doteq\overline{(\frac{1}{z}) \circ \pi \circ \psi  } (x -iy)$ for
$x+iy$ in the lower rectangle, $(\overline{\ \ })$ is the complex conjugation.
As a result $\phi$ is a local biholomorphism in the upper and low open
rectangles and $\phi( x + i0 )\equiv0$, that is a contradiction.

\noindent
Case 2. The path $\beta \in \CC$ is a straight line segment.
There exists a rotation $\e^{i\theta}$, such that $\psi^{-1}(\beta)$
coincides with a trajectory of $\e^{i\theta}X$, see Figure \ref{3-fig}.b.
Hence we can apply the above argument.

\noindent
Case 3. The path $\beta$ has arbitrary shape in $\CC$. Assume by a 
moment that there is a trajectory $z(t)$ of $\e^{i\theta}X$ that 
determine a secant of $\beta$, see Figure \ref{3-fig}.c. By case 2
the conformal structure on the surface $(M'\subset \RR^2,J)$ bounded by
$z(t)$ has conformal hole. The region $(U'J)$ bounded by $z(t)$
and $\beta$ is biholomorphic to the Poincar\'e disk $\Delta$.
If we consider the paste of $(M'\cup \Delta, J)$ obviously also has
a conformal hole.

The finite trajectory gap concept can be used for vector field germs $X = \Re{\XX}$ on 
$((\RR^2\backslash \{\overline{0}\},\overline{0}), J)$, where the ideal boundary
to be consider is $\overline{0}$. In this case, the existence of
a finite trajectory gap, means that $0$ is a conformal hole of 
$((\RR^2\backslash \{\overline{0}\},\overline{0}), J)$. 
\hfill $\Box$

\section{Complex structures around isolated singularities} 

Recall assertion (1) in Corollary \ref{teorema-3-introduccion-caracterizacion-GFB}.
Let $H \subset \RR^2 $ be a topological hyperbolic sector 
of $X$ above, 
having as boundary a vertex at 
$p \in \mathfrak{P} \cup \{\infty\}$ and  
two separatrix trajectories 
$\zeta_1$, $\zeta_2 \subset \partial H \cap \RR^2$, 
with $\alpha$ and $\omega$--limits at $p$ 
(or viceversa).
Let $q_j \in \zeta_j$ 
$j=1, 2$, be two nonsingular points of $X$ and consider $C^\infty$ 
embedded transversals to $X$, 

\centerline{$
\sigma_j:[0, \epsilon) \longrightarrow \Sigma_j \subset H , \ \
\sigma_j (0) = q_j
$}

\noindent 
as manifolds with boundary.
There exists a {\it holonomy map} of $X$ 
\begin{equation}\label{holonomia}
hol: (\Sigma_1, q_1) \longrightarrow (\Sigma_2, q_2), \ \
(x,y) \longmapsto \phi_{\tau(x,y)}(x,y),\\
\end{equation}

\noindent 
here  $\phi_\tau(\ , \ )$ is the flow of $X$ and
$\tau = \tau(x,y)$ is a suitable time function on 
$\Sigma_1 \backslash \{ q_1 \}$.
The map $hol$ is a $C^\infty$ diffeomorphism 
in the interior points of $\Sigma_1, \Sigma_2$.

\begin{figure}[h!]
\begin{center}
\scalebox{0.7}{\includegraphics{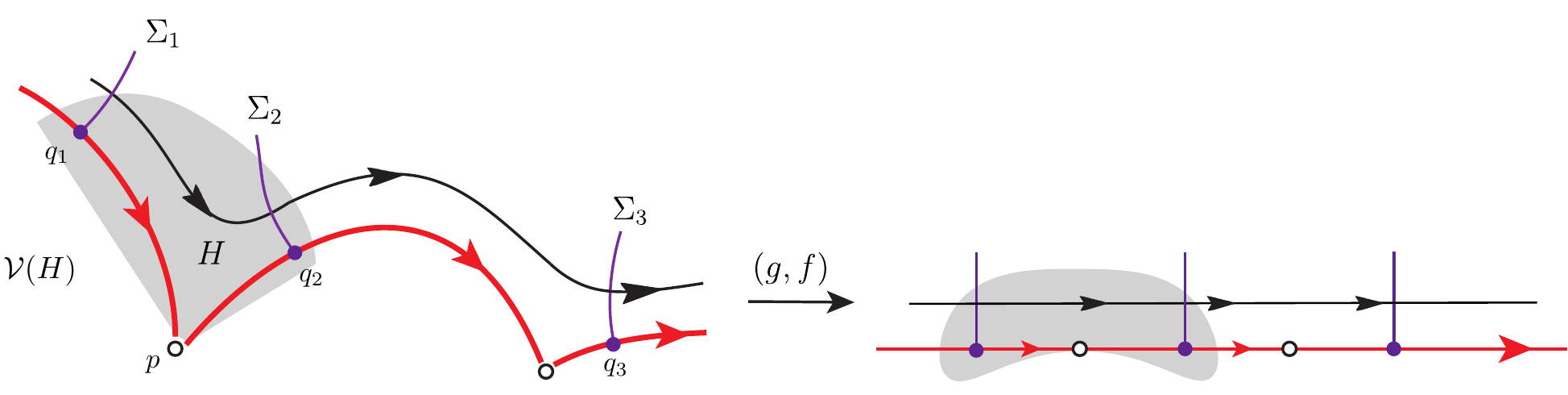}}
\caption{The composition of holonomy maps, must be 
a $C^\infty $ diffeomorphism between transversals 
including their boundary points $q_1, q_2 \in \Gamma(X)$ 
inside of separatrices arriving or leaving 
$p \in \mathfrak{P} \cup \{\infty\}$, according to
Definition \ref{holonomia-difeomorfismo}.}
\label{maquina}
\end{center}
\end{figure}

\begin{definition}\label{holonomia-difeomorfismo}
\begin{upshape}
The {\it holonomy of an hyperbolic sector $H$ 
is a $C^\infty$ diffeomorphism} when the germ 
$hol:(\Sigma_1, q_1) \longrightarrow (\Sigma_2, q_2)$ 
is a $C^\infty$ diffeomorphism of manifolds 
with boundary, see Figure \ref{maquina}.
\end{upshape} 
\end{definition}

\noindent
As far as we known, Definition  
\ref{holonomia-difeomorfismo}
is due to W.\,Kaplan \cite{Kaplan4} p.\,224, 
it was called {\it evenly spread}, 
in J.\,L.\,Weiner \cite{Weiner} p.\,201 
is {\it $\infty$--normal}, 
and also appeared in H.\,Zoladek \cite{Zoladek}.

Corollary 1 is now clear, we provide 
some interesting examples related to real vector fields.

\begin{example} A prototype.
\label{holonomia-hiperbolico-lineal} 
\upshape 
The holonomy of a linear hyperbolic sector $H$ of 

\centerline{ 
$ X (x,y) = \lambda_1 x \del{}{x} + \lambda _2 y\del{}{y}$,
$\lambda_1 \lambda_2 < 0$,
}

\noindent  
is a $C^\infty$ diffeomorphism if and only if 
$\lambda_1 = - \lambda_2$. 
Moreover, in this case the first integral is $f(x,y)= xy$ 
and the infinitesimal symmetry is

\centerline{
$Y (x,y) = \frac{1}{x^2+ y^2}\big(y \del{}{x} + x\del{}{y} \big)$.
}

\end{example}

\begin{example}
A holonomy map which is not a 
$C^\infty$ diffeomorphism,  M.--P. Muller \cite{Muller1}.  
\begin{upshape}
Let us consider

\centerline{
$
X(x,y)= -x^4 \del{}{x} + (x^3 y + 2xy +2) \del{}{y} \  \in C^\infty (\RR^2).
$
}

\noindent 
The holonomy of the hyperbolic sector $H$, 
with vertex in $\infty  \in \SS^2$
which is bounded by the separatrices 
$
\zeta_1 = \{ (0, \tau ) \ \vert \ \tau <0\}$, 
$
\zeta_2 = \{ (\tau, -1/ \tau) \ \vert \ 0< \tau <1  \ \}, 
$
of the vector field $X$ is not a diffeomorphism. 
For this computation Muller used the Liouvillian 
first integral 
$f(x,y)= (xy+1) \e^{- 1/ x^2}$.
The infinitesimal symmetry 
$Y=\nabla {f} / \Vert \nabla {f}  \Vert^2$ 
is not $C^1$ at $\zeta_1$.
\end{upshape}
\end{example}

\begin{example}  {\it The cusp;
a removable singularity.}
\begin{upshape}
The Hamiltonian vector field 

\centerline{$X_f(x,y) = my^{m-1} \del{}{x} + nx^{n-1} \del{}{y} \ \ \ \hbox{on } \mathbb{R}^2$}

\noindent 
of the function $f(x,y)= x^n - y^m$, 
where $n, m \geq 2$ and $(n,m)=1$, 
have at $\overline{0}$ a 2--saddle. 
The union of the separatrices is the singular cusp 
$\{ x^n-y^n=0 \}$ at $\overline{0}$. 
If $n$ is even   
$\rho(x,y)=\e^{-x}/(mx^{m-1}+ny^n)$ 
is a scaling factor for $X_f$ and there exists 
a second vector field  

\centerline{
$Y_g(x,y) = -\e^x \del{}{x} + y\e^x \del{}{y}$
} 
\noindent
linearly independent with $X_f$, such that $
[\rho X_f, \rho Y_g] \equiv 0$ on $\mathbb{R}^2$.
The Proposition \ref{integrabilidad-caja-flujo-global} 
applies and there exists 
a global flow box 

\centerline{$\Psi(x,y)=(y\e^x, x^m-y^n)$}  

\noindent 
for $\rho X_f$ on $\RR^2 \backslash \{0\}$. 
Analogously, if $n$ is odd, there exist a scaling factor 
$\rho(x,y)=\e^{-y}/(mx^m+ny^{n-1})$ 
and a global flow box 
$\Psi(x,y)=(x\e^{y}, x^m-y^n)$ 
for $\rho X_f$. 
\end{upshape}
\end{example}
\begin{example}
{\it Topological saddles with zero linear part.}
\begin{upshape}
The vector field

\centerline{$X(x,y)=x^3\del{}{x}-y^3\del{}{y}$}

\noindent
determines a topological
saddle on $(\RR^2,\overline{0})$. Its foliation $\mathcal{F}(X)$ is symmetric
respect to reflection on both axes, 
hence the holonomy of each hyperbolic sector 
is a $C^\infty$ diffeomorphism on a punctured
neighbourhood of $\overline{0}$. 
There exists a single valued local diffeomorphism
$
\Xi:(\RR^2\backslash\{\overline{0}\},\overline{0})
\longrightarrow 
(\RR^2-\{\overline{0}\},\overline{0})$
such that
$$
\Xi_*(\rho_2(x,y) X(x, y))= x \del{}{x}-y\del{}{y},
$$
for suitable $\rho_2$. The flow box around $\overline{0}$ exists, and is 
$ \Xi \circ \Psi$. 
\end{upshape}
\end{example}

\begin{example}
\begin{upshape}
{\it The saddle node.} 
Let $X \in \mathfrak{X}^\infty (\RR ^2, \overline{0})$ be a saddle node

\centerline{ $X(x,y)= x^2\del{}{x} - \lambda y \del{}{y}, 
\ \ \ \lambda \in \RR^+$.
}

\noindent
The function ${\tt f}(x,y)= e^{\lambda/x}/ y$
is a Liouvillian first integral. The holomomy of the hyperbolic 
sectors is not a $C^\infty$ diffeomorphism.   
Hence, even in the punctured germ domain 
$(\RR^2 \backslash \{ \overline{0}\}, \overline{0})$ 
does not exist an adapted complex structure making $X$
the real part of a singular complex analytical vector field. 
Using a reparametrization

\centerline{
$
\rho(x,y) X(x, y) =\frac{1}{(x^2+y^2)}
\left(
x^2\del{}{x}-\lambda y\del{}{y} 
\right),
$
}
\noindent 
the trajectories of $\rho X$ arrives $\overline{0}$ at finite time.
Moreover, if we remove the origin and the positive real $x$ axis;
there exists an adapted complex structure making 
to $\rho X$ the real part of a holomorphic vector field on 
$(\RR^2 \backslash \{(x, 0) \ \vert \ x \geq 0 \}, J)$, 
see Figure \ref{silla-nodo-figura}. 
In fact, we
consider an embedded transversal to $\rho X$, 

\centerline{
$
\sigma:(-\epsilon, \epsilon) \longrightarrow \Sigma_1 \subset 
\RR^2 \backslash \{ (x, y) \ \vert \ (x, 0), x \geq 0 \}$.}

\noindent
The vector field $Y$ tangent to $\Sigma_1$ is transversal to $\rho X$.
We extend this transversal data over the whole 
domain $\RR^2 \backslash \{ (x, 0) \ \vert \ 0\leq x  \}$
using the flow of $\rho X$; this produces a vector field $Y$, such that $[\rho X,Y]=0$.
Figure \ref{silla-nodo-figura} shows the target of the global flow box
$\Psi=(g,f)$ and the shape of the associated flat surface. 
When we identify the two right boundaries of the regions $H$, 
the origin 0 becomes a finite trajectory gap (a conformal hole) 
of the resulting $(\RR^2 \backslash \{ 0 \}, J )$ where  
$J(\rho X)= Y$ and $J(Y)=-\rho X$.

\begin{figure}[h!]
\begin{center}
\scalebox{0.4}{\includegraphics{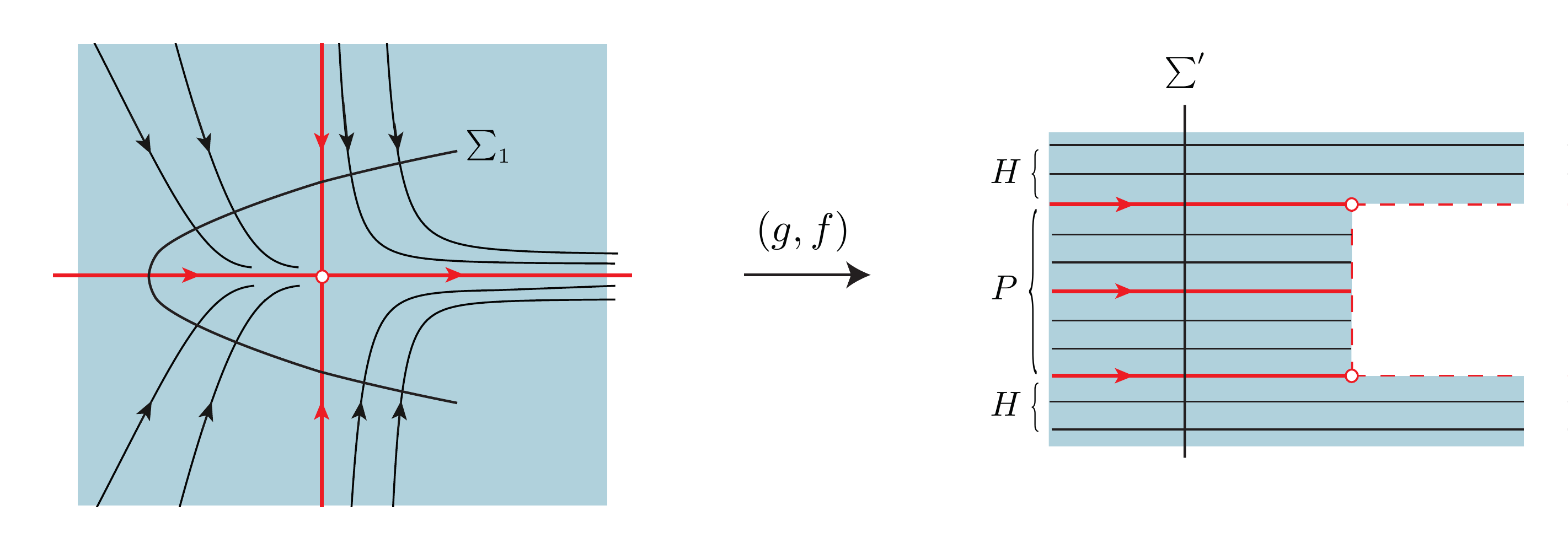}}
\caption{The saddle node $X(x,y)=x^2\del{}{x} -\lambda y \del{}{y}$  determines the
word $HPH$ at $\overline{0}$;
a global flow box exists on the plane minus a
ray $\{ (x, 0) \ \vert \ 0 \leq x \}$, here $\rho X$ becomes the 
real part of a holomorphic vector field.
\label{silla-nodo-figura}
}
\end{center}
\end{figure}

\end{upshape}
\end{example}
The converse of Corollary \ref{teorema-3-introduccion-caracterizacion-GFB}
is the goal of a future work. 


\end{document}